\documentclass[graybox,envcountchap,sectrefs]{svmono}
\usepackage{type1cm}
\usepackage{bpproc}
\usepackage{makeidx}   
\usepackage{graphicx}
\usepackage{multicol}
\usepackage[bottom]{footmisc}
\usepackage{newtxtext}  
\usepackage{newtxmath}

\usepackage{amsmath,hyperref, mdframed}
\usepackage{mathrsfs} 
\usepackage{tikz}
\usepackage[all]{xy}
\usepackage{subfig}
\usepackage{tikz}
\usepackage{subfig}
\usepackage{tabularx}
\usepackage{enumitem}
\usepackage{amscd}
\usepackage{amsfonts}
\usepackage{fancyhdr} 
\usepackage{pgf,tikz}
\usepackage{tikz,tikz-3dplot}
\usepackage{caption}
\usepackage{tikz-cd}
\usepackage{tikzscale}
\usepackage{rotating}
\usepackage{xcolor}
\usepackage{lscape}
\usepackage{xspace}
\usepackage{xfrac}
\usepackage{blindtext}
\usepackage[ruled,vlined]{algorithm2e}
\usepackage[export]{adjustbox}
\usepackage[dvipsnames]{pstricks}
\usepackage{pst-solides3d}
\usepackage{color}
\usepackage{hyperref}
\hypersetup{
    colorlinks=true, 
    linktoc=all,     
    linkcolor=blue,  
    citecolor=blue,
    filecolor=blue,
    urlcolor=blue
}

\usetikzlibrary{arrows}
\usetikzlibrary{decorations.markings}
\usetikzlibrary{decorations.markings}
\usetikzlibrary{backgrounds,shapes}
\usetikzlibrary{decorations.markings}
\usetikzlibrary{patterns}
\usetikzlibrary{calc,perspective}

    \pgfdeclarelayer{background} 
    \pgfdeclarelayer{foreground}
    \pgfsetlayers{background,main,foreground}

\definecolor{aliceblue}{rgb}{0.9, 0.95, 1.0}

\definecolor{violet1}{RGB}{69, 14, 71}
\definecolor{color2}{RGB}{162, 165, 149}
\definecolor{color3}{RGB}{180, 162, 138}

\definecolor{cof}{RGB}{219,144,71}
\definecolor{pur}{RGB}{186,146,162}
\definecolor{greeo}{RGB}{91,173,69}
\definecolor{greet}{RGB}{52,111,72}

\newcommand{\mytriangle}[4] 
{
  \coordinate (center) at ($1/3*(#1)+1/3*(#2)+1/3*(#3)$);
  \coordinate (m12)    at ($(#1)!0.5!(#2)$);
  \coordinate (m13)    at ($(#1)!0.5!(#3)$);
  \coordinate (m23)    at ($(#2)!0.5!(#3)$);
  \draw[fill=pink] (center) -- (m12) -- (#1) -- (m13) -- cycle;
  \draw[fill=red] (center) -- (m12) -- (#2) -- (m23) -- cycle;
  \draw[fill=purple] (center) -- (m13) -- (#3) -- (m23) -- cycle;
  \draw[thick,fill=black,fill opacity=#4]  (#1) -- (#2)  -- (#3) -- cycle;
}

\numberwithin{equation}{section}

\makeatletter
\def\namedlabel#1#2{\begingroup
   \def\@currentlabel{#2}%
   \label{#1}\endgroup
}
\makeatother

\makeindex

\begin{document}

\theoremstyle{plain}                    
\newtheorem*{khthm}{Killing-Hopf's Theorem}
\newtheorem*{pthm}{Poincar\'e's Theorem}
\newtheorem*{ethm}{Euler's observation}
\newtheorem{thma}{Theorem}
\newtheorem{propa}[thma]{Proposition}
\newtheorem{cora}[thma]{Corollary}

\renewcommand{\thethma}{\Alph{thma}}
\newtheorem{thmaa}{Theorem}
\renewcommand{\thethmaa}{\Alph{thmaa}.2}
\newtheorem*{thmnn}{Theorem}
\newtheorem*{propnn}{Proposition \ref{aiffd}}
\newtheorem{conj}{Conjecture}	

\theoremstyle{definition}
\newtheorem*{defnn}{Definition}
\newtheorem*{qs}{Question}
\newcommand{\bluecomment}[1]{\textcolor{blue}{#1}}
\newcommand*{\bigchi}{\mbox{\large$\chi$}}

\graphicspath{{./}{images/}}

\newcommand{\mf}{\mathfrak}
\newcommand{\spec}{\mathrm{Spec}}

\title{Hilbert modular forms and Galois representations}
\author{Ajith Nair\\ 
Ajmain Yamin}
\thispagestyle{empty}

\maketitle
\tableofcontents

\chapter{Hilbert modular forms \& Galois representations}

\section{Introduction}

The notion of a Hilbert modular form originates from the unpublished work of David Hilbert and the  habilitation thesis of his doctoral student, Otto Blumenthal \cite{blumenthal:1903}. The theory gained prominence after development of the theory of several complex variables. Shimura laid the foundation for the arithmetic theory of Hilbert modular forms \cite{shimura:1978}, and Hida provided a framework for the Hecke theory \cite{hida:1988}.

The Taniyama-Shimura-Weil conjecture states that every elliptic curve over $\Q$ is modular. The conjecture was proved by Wiles for semistable elliptic curves, thereby proving Fermat's last theorem, and the full theorem was proved by Breuil, Conrad, Diamond and Taylor.

Motivated by the success of Wiles' proof, Darmon \cite{darmon:2000} initiated a program to tackle the Diophantine equation
\begin{equation}
\label{eq:GFE}
    x^p + y^q = z^r
\end{equation}
where $p$, $q$ and $r$ are primes and $p$ is odd. Darmon's approach begins by associating certain Frey representations to putative solutions to Equation \ref{eq:GFE} which correspond to certain abelian varieties of $GL(2)$-type over totally real fields.  The next step in Darmon's program is to establish the modularity of these abelian varieties.  Here, \emph{modularity} refers to the existence of an appropriate Hilbert modular form whose associated Galois representation corresponds to the abelian variety.  The following conjecture generalizes the Taniyama-Shimura-Weil conjecture, and asserts that abelian varieties of $GL(2)$-type over totally real fields are modular.

\begin{conjecture}\label{conj:modularityGL2AbVar}
    If $A$ is an abelian variety of $GL(2,E)$-type over a number field $K$ of conductor $\mathfrak{n}$, then there exists a Hilbert modular form $f$ over $K$ of weight $2$ and level $\mathfrak{n}$ such that
    $$   \rho_{f,\lambda} \cong \rho_{A,\lambda} $$
    for all primes $\lambda$ of $E$.
\end{conjecture}
Here $\rho_{f,\lambda}$ (respectively, $\rho_{A,\lambda}$) are compatible systems of $\lambda$-adic representations of $K$ associated to the Hilbert modular form $f$ (respectively, abelian variety $A$).

The conjecture has been proved for elliptic curves over certain classes of number fields. For more details, we refer the reader to the article \emph{Elliptic curves over totally real fields: A survey} by Bidisha Roy and Lalit Vaishya in this proceedings. The recent work of Billerey, et al. \cite{billerey:2023} implements and extends aspects of Darmon's program to solve several cases of the the generalized Fermat equation \ref{eq:GFE}.

The purpose of this article is to introduce the theory of Hilbert modular forms and Galois representations, and describe what it means to attach a compatible system of Galois representations $\rho_{f,\lambda}$ to a Hilbert modular form $f$. The organization of this article is as follows. In Section \ref{section:HMFclassical}, we discuss the theory of Hilbert modular forms in the classical language of holomorphic functions on $n$-copies of the hyperbolic plane. In Section \ref{section:HMFadelic}, we re-introduce Hilbert modular forms in the adelic formalism, allowing for a smooth treatment of their Hecke theory and $L$-functions. In Section \ref{section:Galoisreps}, we give an introduction to Galois representations and their $L$-functions. We conclude this article with a discussion of the landmark result, Theorem \ref{theorem:mainresult}, which asserts the existence of a compatible system of Galois representations associated to a Hilbert modular form.

\section{Hilbert modular forms (classical formulation)}
\label{section:HMFclassical}

This section is adapted from Chapter 2 of \emph{The 1-2-3 of Modular Forms} by Jan Hendrik Bruinier \cite{zagier123:2008} and Chapter 1 of Paul Garrett's \emph{Holomorphic Hilbert Modular Forms} \cite{garrett:1989}.

\subsection{Hilbert modular group}\label{section:hilbertModularGroup}
Let $F$ be a totally real field of degree $n$ over $\mathbb{Q}$.  Write $\mathcal{O}_F$ for the ring of integers in $F$.  The \emph{Hilbert modular group} $$\Gamma_F := SL(2,\mathcal{O}_F)$$ embeds into $SL(2,\mathbb{R})^n$ via the $n$ real embeddings $F \hookrightarrow \mathbb{R}$, denoted $x \mapsto x^{(j)}$ for $j=1,\cdots, n$.  There is an action $\Gamma_F \curvearrowright \mathbb{H}^n$ of the Hilbert modular group on $n$ copies of the hyperbolic plane $\mathbb{H} := \{z\in \mathbb{C} \mid \Im(z) > 0 \}$ given by embedding $\Gamma_F$ into $SL(2,\mathbb{R})^n$ and letting each copy of $SL(2,\mathbb{R})$ act on the corresponding copy of $\mathbb{H}$ via linear fractional transformations $z\mapsto (az+b)/(cz+d)$.

In this section, we will be concerned with subgroups of $SL(2,F)$  which are \emph{commensurable} with the Hilbert modular group, i.e., subgroups $\Gamma \subset SL(2,F)$ such that $\Gamma \cap \Gamma_F$ has finite index in both $\Gamma$ and $\Gamma_F$.  For each ideal $\mathfrak{n}$ of $\mathcal{O}_F$, the \emph{principal congruence subgroup of level $\mathfrak{n}$} $$\Gamma_F(\mathfrak{n}) := \ker(\Gamma_F \rightarrow SL(2,\mathcal{O}_F/\mathfrak{n}))$$ is an example of a commensurable subgroup.  A \emph{congruence subgroup} is a subgroup of $\Gamma_F$ which contains a principal congruence subgroup.  All finite index subgroups of the Hilbert modular group of a totally real number field of degree greater than 1 are congruence subgroups \cite{serre:1970}.

\subsection{Hilbert modular varieties}
Throughout this section, $\Gamma$ will denote a congruence subgroup of $\Gamma_F$.  The group $\Gamma$ acts discretely on $\mathbb{H}^n$. The quotient $$Y(\Gamma) := \Gamma \backslash \mathbb{H}^n$$ is a non-compact Hausdorff topological space, which we refer to as a \emph{Hilbert modular variety}.  If the image of $\Gamma$ in $PSL(2,F)$ is torsion free, then the quotient map $\mathbb{H}^n \rightarrow Y(\Gamma)$ is a covering projection and the Hilbert modular variety $Y(\Gamma)$ is naturally an $n$-dimensional complex manifold.  In general, $\mathbb{H}^n \rightarrow Y(\Gamma)$ is a branched covering projection, with finitely many branch points in $Y(\Gamma)$.  

Since the measure $\prod_{j=1}^n\frac{dx_j dy_j}{y_j^2}$ is invariant under the action $SL(2,\mathbb{R})^n \curvearrowright \mathbb{H}^n$ it descends to a well-defined measure on $Y(\Gamma)$.  It can be shown that the $Y(\Gamma)$ has finite volume with respect to this measure \cite[pg. 46]{garrett:1989}.

Hilbert modular varieties are non-compact but can be compactified by adding a finite number of cusps.   A \emph{cusp} of $\Gamma$ is an equivalence class of $\mathbb{P}^1(F):= F \cup \{\infty\}$ under the action of $\Gamma$.  The cusps of $\Gamma_F$ are in one to one correspondence with ideal classes of $F$.  Indeed, the map \begin{align*}
\varphi: \Gamma_F \backslash \mathbb{P}^1(F) &\rightarrow \operatorname{Cl}(F)\\
[\alpha : \beta] &\mapsto \alpha \mathcal{O}_F + \beta \mathcal{O}_F
\end{align*}
is a bijection \cite[pg. 7]{garrett:1989}. Consequently, any subgroup $\Gamma \subset SL(2,F)$ commensurable with the Hilbert modular group has finitely many cusps. 

\subsection{Bailey-Borel compactification}
Consider $$(\mathbb{H}^n)^* := \mathbb{H}^n \cup \mathbb{P}^1(F) = \mathbb{H}^n \cup F \cup \{\infty \}.$$  It can be equipped with an appropriate topology which extends that of $\mathbb{H}^n$.  We describe this topology by specifying a neighborhood basis for each point of $\mathbb{P}^1(F)$. 
A neighborhood basis for the point at infinity is given by the sets $$N_C :=\{\infty\} \cup \{ z \in \mathbb{H}^n \mid {\textstyle\prod_{j=1}^n}\Im(z_j) > C\}$$ as $C>0$ ranges over all positive reals.  A neighborhood basis for a point $\alpha \in \mathbb{P}^1(F)$ is given by $\{\gamma(N_C) \}_{C>0}$ where $\gamma \in \Gamma$ is any group element such that $\gamma(\infty) = \alpha$. Such an element $\gamma$ always exists as $SL(2,F)$ acts transitively on $\mathbb{P}^1(F)$ via linear fractional transformations $[\alpha : \beta] \mapsto [a\alpha + b\beta : c\alpha + d\beta]$.

The \emph{Baily-Borel compactification} $X(\Gamma)$ of the Hilbert modular variety $Y(\Gamma)$ is the quotient space defined as $$X(\Gamma) := \Gamma \backslash (\mathbb{H}^n)^*.$$  Here the action $\Gamma \curvearrowright (\mathbb{H}^n)^*$ is given by letting $\Gamma$ act on $\mathbb{H}^n$ and $\mathbb{P}^1(F)$ separately, as described above.  The Hilbert modular variety $Y(\Gamma)$ is included in its Baily-Borel compactification $X(\Gamma)$ with complement being the finite set of cusps of $\Gamma$.  In contrast to $Y(\Gamma)$ being non-compact, the space $X(\Gamma)$ is a compact Hausdorff topological space.  When $n>1$, even if $Y(\Gamma)$ is a complex manifold, the space $X(\Gamma)$ can never be given the structure of a complex manifold.  This is because no neighborhood of a cusp in $X(\Gamma)$ is homeomorphic to an open ball in $\mathbb{R}^{2n}$ when $n>1$ \cite[pg. 30]{freitag:1990}. Instead, $X(\Gamma)$ can be can be given the structure of a normal complex space and also, remarkably, the structure of a projective algebraic variety \cite[Ch. 2, Sect. 4]{freitag:1990}.  Consequently, Hilbert modular varieties can be seen as quasi-projective varieties.

\subsection{Example: $\Gamma = \Gamma_F[\mathfrak{n}]$, where $F = \mathbb{Q}(\sqrt{5})$ and $\mathfrak{n}=2\mathcal{O}_F$}
When $n=2$, we refer to Hilbert modular varieties as \emph{Hilbert modular surfaces}.

Consider the case when $F = \mathbb{Q}(\sqrt{5})$ is the real quadratic field of smallest positive discriminant.  The principal congruence subgroup of level two, $\Gamma = \Gamma_F[\mathfrak{n}]$, $\mathfrak{n} = 2 \mathcal{O}_F$, has torsion-free image in $PSL(2,F)$.  Thus the Hilbert modular surface $Y(\Gamma)$ is a complex surface.  Note $F$ has class number one and $\Gamma$ has index four in $\Gamma_F$.  The group $\Gamma$ has five cusps, represented by $0,1,\infty, \frac{1+\sqrt{5}}{2},\frac{1-\sqrt{5}}{2}$.  The Bailey-Borel compactification $X(\Gamma)$ is no longer a complex manifold, but can be embedded into complex projective space as a projective algebraic variety.  Indeed, consider the \emph{(weight two) Eisenstein series} defined by 
\begin{equation}\label{eq:wt2eisenstein}
E_0(z_1,z_2) :=
\lim_{s \searrow 0} 
\sum_{\substack{c,d \in \mathcal{O}_F \\
c\mathcal{O}_F + d\mathcal{O}_F= \mathcal{O}_F\\ c\equiv 0, d\equiv 1 \mod \mathfrak{n}}}' \frac{\left|c^{(1)}z_1 + d^{(1)}\right|^{-s} \left|c^{(2)}z_2 + d^{(2)}\right|^{-s}}{(c^{(1)} z_1 + d^{(1)})^2 (c^{(2)} z_2 + d^{(2)})^2}
\end{equation}
where the dash indicates that the summation is modulo the group $\mathcal{O}_F^{\times 3}$ of units $\varepsilon \equiv 1 \mod \mathfrak{n}$.  The Eisenstein series $E_0(z_1,z_2)$ converges and vanishes at all cusps of $\Gamma$ except for $\infty$ \cite[pg. 173]{freitag:1990}.  There are four other Eisenstein series $E_1, E_2, E_3, E_4$ associated to the other four cusps in a similar way.  The mapping $X(\Gamma) \rightarrow \mathbb{P}^4(C)$ given by $$(z_1,z_2) \mapsto [E_0(z_1,z_2) : E_1(z_1,z_2) : E_2(z_1,z_2) : E_3(z_1,z_2) : E_4(z_1,z_2)]$$ is well-defined and embeds $X(\Gamma)$ into $\mathbb{P}^4(\mathbb{C})$ as \emph{Klein's icosahedral surface}, the projective algebraic variety in $\mathbb{P}^4(\mathbb{C})$ defined by the equations $$\sigma_2(u_0,u_1,u_2,u_3,u_4) = \sigma_4(u_0,u_1,u_2,u_3,u_4)=0$$ where $\sigma_i$ is the $i$'th symmetric polynomial in five variables \cite{hirzebruch:1976}. 

\subsection{Hilbert modular forms}
\label{section:Hilbertmodularforms}
Let $\mathbf{k} = (k_1,\cdots, k_n)\in \mathbb{Z}^n$.  For any function $f:\mathbb{H}^n \rightarrow \mathbb{C}$ and any $\gamma = \begin{pmatrix}
a & b \\
c & d 
\end{pmatrix}\in SL(2,F)$, denote by $f\vert_\mathbf{k} \gamma$ the function $$f\vert_\mathbf{k} \gamma : \mathbb{H}^n \rightarrow \mathbb{C} \quad \text{given by} \quad z \mapsto \prod_{j=1}^n(c^{(j)}z_j + d^{(j)})^{-k_j}f(\gamma z).$$

For $n>1$, a \emph{Hilbert modular form} for $\Gamma$ of weight $\mathbf{k}$ is a holomorphic function $f: \mathbb{H}^n \rightarrow \mathbb{C}$ which obeys the transformation law $f\vert_\mathbf{k} \gamma = f$ for each $\gamma \in \Gamma$.  In particular, if $\Gamma = \Gamma_F$ then $f(z + a) = f(z)$ for each $a \in \mathcal{O}_F$. The \emph{translation module} $\mathfrak{t} = \mathfrak{t}_\Gamma$ for $\Gamma$ is the set of translations in $\Gamma$, i.e., $$\mathfrak{t}_\Gamma := \{ a\in \mathbb{R}^n \mid \exists \gamma \in \Gamma \text{ such that } \gamma(z) = z + a \}.$$ It is a lattice of rank $n$.  It is clear that $f(z + a) = f(z)$ for all $a\in \mathfrak{t}_\Gamma$.  This fact allows us to obtain a \emph{Fourier expansion} for $f$, viz.,
\begin{align}
\label{eq:fourierSeries}
f(z) &= \sum_{\nu \in \mathfrak{t}^\vee} A(\nu) e^{2\pi i \langle \nu , z \rangle}, \\
\label{eq:fourierCoefficient}
A(\nu) &:= \frac{1}{\operatorname{covol}(\mathfrak{t})}\int_{\mathbb{R}^n/ \mathfrak{t}} f(x + iy) e^{- 2\pi i \langle \nu, x + i y \rangle} \prod_{j=1}^n dx_j,
\end{align}
where $\langle \cdot, \cdot \rangle$ is the usual inner product on $\mathbb{R}^n$ and $\mathfrak{t}^\vee$ is the \emph{dual lattice} $$\mathfrak{t}^\vee := \{\nu \in \mathbb{R}^n \mid \langle \nu, a \rangle \in \mathbb{Z} \text{ for all } a\in \mathfrak{t}\}.$$  The number $A(\nu)$ is called the $\nu$'th \emph{Fourier coefficient} for $f$.  The fact that the above formula for $A(\nu)$, Equation \ref{eq:fourierCoefficient}, is independent of the choice of $y\in \mathbb{R}_+^n$ follows from the holomorphicity of $f$.  The \emph{Koecher principle} states that $A(\nu) =0 $ unless $\nu$ is totally positive or zero \cite[pg. 10]{garrett:1989}. 

The above expansion, Equation \ref{eq:fourierSeries}, is the Fourier expansion of $f$ at $\infty$.  We may also obtain a Fourier expansion for $f$ at other cusps of $\Gamma$.  Indeed, if $\kappa \in F$ is another cusp of $\Gamma$, there exists some $\gamma \in SL(2,F)$ for which $\gamma \kappa = \infty$.  It may be seen that $\gamma \Gamma \gamma^{-1}$ is still commensurable with the Hilbert modular group and that $f\vert_{\mathbf{k}}\gamma$ is a Hilbert modular form for $\gamma \Gamma \gamma^{-1}$.  The Fourier expansion of $f\vert_{\mathbf{k}}\gamma$ at $\infty$ will be considered as a Fourier expansion of $f$ at $\kappa$, though it is not unique in that it depends on the choice of $\gamma$.  Nevertheless, whether or not the constant term in this Fourier expansion is zero is independent of the choice of $\gamma$.  We say $f$ \emph{vanishes at the cusp} $\kappa$ if the constant term of any Fourier expansion at $\kappa$ is $0$.  We say $f$ is a \emph{cusp form} for $\Gamma$ if $f$ vanishes at all cusps of $\Gamma$.

\subsection{Finite dimensionality}
Let $\mathbf{k}=(k_1,\cdots,k_n)\in \mathbb{Z}^n$ be a vector of weights.  We write $M_\mathbf{k}(\Gamma)$ to denote the vector space of Hilbert modular forms for $\Gamma$ of weight $\mathbf{k}$, and $S_\mathbf{k}(\Gamma)$ to denote the subspace of cusp forms.  The codimension of $S_\mathbf{k}(\Gamma)$ in $M_\mathbf{k}(\Gamma)$ is bounded by the number of cusps of $\Gamma$. Indeed, if $\kappa_1, \cdots, \kappa_h$ are the distinct cusps of $\Gamma$, and $\gamma_1, \cdots, \gamma_h \in SL(2,F)$ are fixed elements such that $\gamma_j \kappa_j = \infty$, we may define a linear map $M_\mathbf{k}(\Gamma) \rightarrow \mathbb{C}^h$ with kernel $S_\mathbf{k}(\Gamma)$ given by $f\mapsto (c_{\gamma_j}(f))_{j=1}^h$ where $c_{\gamma_j}(f)$ is the constant term in the Fourier expansion of $f\vert_{\mathbf{k}} \gamma_j$ at $\infty$.  Thus, to show finite dimensionality of these vector spaces, it suffices to show $S_\mathbf{k}(\Gamma)$ is finite dimensional.  It is also true that when $\mathbf{k}$ is not parallel, i.e. $\mathbf{k} \neq k (1,\cdots, 1)$ for some $k\in \mathbb{Z}$, then $M_\mathbf{k}(\Gamma) = S_\mathbf{k}(\Gamma)$.  This is a consequence of the Koecher principle \cite[pg. 10]{garrett:1989}.  Moreover, if $\mathbf{k}\neq 0$ and any of its entries are non-positive, then $M_{\mathbf{k}}(\Gamma) = \{ 0\}$ \cite[pg. 18]{geer:1987}.  Note, $M_0(\Gamma) = \mathbf{C}$ consists of constant functions.

For any $f,g \in S_\mathbf{k}(\Gamma)$, the \emph{Petersson inner product} $\langle f, g \rangle$ is defined as
\begin{equation}
\langle f,g \rangle := \int_{\Gamma \backslash \mathbb{H}^n} f(z)\overline{g(z)}\prod_{j=1}^n y_j^{k_j} \frac{dx_jdy_j}{y_j^2}.
\end{equation}
This defines an inner product on $S_\mathbf{k}(\Gamma)$ and makes $S_\mathbf{k}(\Gamma)$ into a Hilbert space.  The finite dimensionality of $S_\mathbf{k}(\Gamma)$ follows from existence of a uniform constant $A$ such that $||f||_\infty \leq A ||f||_2$, where $||f||_\infty$ is the $L^\infty$-norm on $S_\mathbf{k}(\Gamma)$ given by $||f||_\infty := \max_{z\in \mathbb{H}^n} |f(z)| \prod_{j=1}^n y^{k_j}$, and $||f||_2$ is the $L^2$-norm on $S_\mathbf{k}(\Gamma)$ coming from the Petersson inner product \cite[pg. 68]{freitag:1990}.  Indeed, making the constant $A$ effective would yield the upper bound on dimension $\dim_{\mathbb{C}}S_{\mathbf{k}}(\Gamma) \leq A^2 \cdot \operatorname{vol}(\Gamma \backslash \mathbb{H}^n)$.

\subsection{Eisenstein series}
For any $2\times 2$ matrix $\gamma$, write $a_\gamma,b_\gamma,c_\gamma, d_\gamma$ to denote the entries of $\gamma$ so that $$\gamma = \begin{pmatrix}
    a_\gamma & b_\gamma \\
    c_\gamma & d_\gamma
\end{pmatrix}.$$

Let $\mathbf{k} = k(1,\cdots,1)$ where $2<k\in \mathbf{Z}$ is even.  Let $\Gamma \subset \Gamma_F$ be a congruence subgroup.  As stated in Section \ref{section:hilbertModularGroup}, when $n\geq 2$ this is equivalent to being of finite index in $\Gamma_F$.  We generally follow Paul Garrett's \emph{Holomorphic Hilbert Modular Forms} \cite{garrett:1989} in defining Eisenstein series for $\Gamma$.  For each $\alpha \in SL(2,F)$ define an \emph{Eisenstein series of weight $\mathbf{k}$} for $\Gamma$ by
\begin{equation}\label{eq:wt>2eisenstein}
E_\alpha(z) := \sum_{\gamma \in S_\alpha} \prod_{j=1}^n(c_{\alpha^{-1}\gamma}^{(j)}z_j + d_{\alpha^{-1}\gamma}^{(j)})^{-k}
\end{equation} 
where $S_\alpha$ is a set of representatives for $(\alpha P \alpha^{-1} \cap \Gamma)\backslash \Gamma$. Here, $P$ denotes the group of upper triangular matrices in $SL(2,F)$.  The series defining $E_\alpha$ does not depend on the choice of representatives $S_\alpha$.  Indeed, if $\gamma$ and $\gamma'$ represent the same left coset in $(\alpha P \alpha^{-1} \cap \Gamma) \backslash \Gamma$, so that $\gamma' = \alpha p \alpha^{-1} \gamma$ for some $p \in P$ satisfying $\alpha p \alpha^{-1} \in \Gamma$, then 
\begin{align*}
\prod_{j=1}^n(c_{\alpha^{-1}\gamma'}^{(j)}z_j + d_{\alpha^{-1}\gamma'}^{(j)})^{-k} &= \prod_{j=1}^m (d_p^{(j)}c_{\alpha^{-1}\gamma}^{(j)}z_j + d_p^{(j)}d_{\alpha^{-1}\gamma}^{(j)})^{-k} \\
&= \mathrm{Nm}^F_\mathbb{Q}(d_p)^{-k} \cdot \prod_{j=1}^n(c_{\alpha^{-1}\gamma}^{(j)}z_j + d_{\alpha^{-1}\gamma}^{(j)})^{-k}, 
\end{align*}
and moreover $\mathrm{Nm}^F_\mathbb{Q}(d_p)^{-k} = 1$ as $d_p$--being an eigenvalue of a matrix in $\Gamma_F$--is a unit in $\mathcal{O}_F$ and $k$ is even.

Recall, a series is \emph{absolutely convergent} if the sum of absolute values of its terms is convergent.  The Eisenstein series defined in Equation \ref{eq:wt>2eisenstein} are absolutely convergent, since $k>2$ \cite[pg. 34]{garrett:1989}.  This allows us to freely rearrange the terms in the series defining the Eisenstein series and observe that $E_\alpha\vert_\mathbf{k} \delta =E_\alpha$ for any $\delta \in \Gamma$.  Indeed, if $\delta = \begin{pmatrix}
e & f \\
g & h 
\end{pmatrix} \in \Gamma$, then right multiplication by $\delta$ induces a permutation on the right cosets in $(\alpha P \alpha^{-1} \cap \Gamma)\backslash \Gamma$, so that, if $S_\alpha$ is a set of coset representatives for $(\alpha P \alpha^{-1} \cap \Gamma)\backslash \Gamma$, then so is $S_\alpha \delta := \{ \gamma \delta \mid \gamma \in S_\alpha\}$.  Hence, 
\begin{align*}
(E_\alpha \vert_{\mathbf{k}}\delta)(z) &= \prod_{j=1}^n (g^{(j)}z_j + h^{(j)})^{-k} \sum_{\gamma} \prod_{j=1}^n\left(c_{\alpha^{-1}\gamma}^{(j)}\left(\frac{e^{(j)}z_j + f^{(j)}}{g^{(j)}z_j + h^{(j)}}\right) + d_{\alpha^{-1}\gamma}^{(j)}\right)^{-k}\\
&= \sum_{\gamma\in S_\alpha} \prod_{j=1}^n ((c_{\alpha^{-1}\gamma}^{(j)}e^{(j)} + d_{\alpha^{-1}\gamma}^{(j)}g^{(j)})z_j + c_{\alpha^{-1}\gamma}^{(j)}f^{(j)} + d_{\alpha^{-1}\gamma}^{(j)}h^{(j)})^{-k}\\
&= \sum_{\gamma \in S_\alpha} \prod_{j=1}^\infty (c_{\alpha^{-1}\gamma \delta}^{(j)}z_j + d_{\alpha^{-1}\gamma\delta}^{(j)})^{-k} \quad \text{where } \begin{pmatrix}
a_{\alpha^{-1}\gamma\delta} & b_{\alpha^{-1}\gamma\delta} \\
c_{\alpha^{-1}\gamma\delta} & d_{\alpha^{-1}\gamma\delta} 
\end{pmatrix} = \alpha^{-1}\gamma\delta \\
&= \sum_{\gamma' \in S_\alpha\delta}\prod_{j=1}^\infty (c_{\alpha^{-1}\gamma'}^{(j)}z_j + d_{\alpha^{-1}\gamma'}^{(j)})^{-k} = E_\alpha(z).
\end{align*}

Thus, Eisenstein series are examples of Hilbert modular forms, and $E_\alpha \in M_\mathbf{k}(\Gamma)$.  Write $E_{\mathbf{k}}(\Gamma)$ to denote the subspace of $M_\mathbf{k}(\Gamma)$ spanned by Eisenstein series.  Notice that $E_\alpha$ is a constant multiple of $ E_\beta$ if $\Gamma \alpha P = \Gamma \beta P$.  Indeed, taking $S_\alpha$ to be a set of coset representatives for $(\alpha P \alpha^{-1} \cap \Gamma)\backslash \Gamma$, we see that if $\beta = \delta \alpha p$ for some $\delta \in \Gamma$ and $p\in P$ then $\delta S_\alpha \delta^{-1} := \{ \delta \gamma \delta^{-1} \mid \gamma \in S_\alpha \}$ is a set of coset representatives for $(\beta P \beta^{-1} \cap \Gamma)\backslash \Gamma$ and \begin{align*}
E_\beta(z) &= \sum_{\gamma' \in \delta S_\alpha \delta^{-1}} \prod_{j=1}^n(c_{\beta^{-1}\gamma'}^{(j)}z_j + d_{\beta^{-1}\gamma'}^{(j)})^{-k}\\ &= \sum_{\gamma \in S_\alpha} \prod_{j=1}^n(c_{(\delta \alpha p)^{-1}(\delta \gamma \delta^{-1})}^{(j)}z_j + d_{(\delta \alpha p)^{-1}(\delta \gamma \delta^{-1})}^{(j)})^{-k}\\
&= \sum_{\gamma \in S_\alpha} \prod_{j=1}^n(c_{p^{-1}\alpha^{-1}\gamma \delta^{-1}}^{(j)}z_j + d_{p^{-1}\alpha^{-1}\gamma \delta^{-1}}^{(j)})^{-k}\\ &=  \sum_{\gamma' \in S_\alpha\delta^{-1}} \prod_{j=1}^n(c_{p^{-1}\alpha^{-1}\gamma' }^{(j)}z_j + d_{p^{-1}\alpha^{-1}\gamma'}^{(j)})^{-k}\\
&= \mathrm{Nm}^F_\mathbb{Q}(d_{p^{-1}})^{-k} \cdot \sum_{\gamma' \in S_\alpha\delta^{-1}} \prod_{j=1}^n(c_{\alpha^{-1}\gamma' }^{(j)}z_j + d_{\alpha^{-1}\gamma'}^{(j)})^{-k}\\ &= \mathrm{Nm}^F_\mathbb{Q}(d_{p^{-1}})^{-k} \cdot E_\alpha(z).
\end{align*}

Let $\{\delta_i\}$ be a set of representatives for $\Gamma \backslash SL(2,F) / P$.  It can be shown that $\{\kappa_i := \delta_i \infty\}$ is a complete set of representatives for the cusps of $\Gamma$ \cite[pg. 8]{garrett:1989}.  Moreover, it can be shown that the Eisenstein series $E_\alpha$ vanishes at all cusps of $\Gamma$ except for the unique cusp corresponding to the double coset $\Gamma \alpha P$ where it does not vanish \cite[pg. 34]{garrett:1989}. It follows that the Eisenstein series $\{ E_{\delta_i}\}$ form a basis for $E_\mathbf{k}(\Gamma)$, and $E_\mathbf{k}(\Gamma)$ has dimension equal to the number of cusps of $\Gamma$.  Moreover, given any Hilbert modular form $f\in M_\mathbf{k}(\Gamma)$, one can subtract off a suitable linear combination $E$ of the Eisenstein series $\{ E_{\delta_i}\}$ so that $f-E\in S_\mathbf{k}(\Gamma)$ is a cusp form.  So we see, $M_\mathbf{k}(\Gamma) = E_\mathbf{k}(\Gamma) \oplus S_\mathbf{k}(\Gamma)$.  Finally, this direct sum decomposition is orthogonal with respect to the Petersson inner product \cite[pg. 51]{garrett:1989}.

The Eisenstein series $E_0$ defined in Equation \ref{eq:wt2eisenstein} is also a Hilbert modular form.  It can be shown that $E_0 \in S_{(2,2)}(\Gamma)$ for $\Gamma = \Gamma_F(\mathfrak{n})$ where $F = \mathbb{Q}(\sqrt{5})$ and $\mathfrak{n} = 2\mathcal{O}_F$ \cite[pg.65 \& pg. 160]{freitag:1990}.  Note that the formula defining Eisenstein series in Equation \ref{eq:wt>2eisenstein} does not converge absolutely for weight $\mathbf{k} = (2, \cdots, 2)$.  This is why Eisenstein series of weight $\mathbf{k} = (2, \cdots, 2)$ are defined using the so-called ``Hecke trick'' or \emph{Hecke regularization}, which involves analytic continuation in an auxiliary parameter $s\in \mathbb{C}, \Re(s) > 0$, followed by evaluation at $s=0$ \cite[pg. 65 \& Chapter III, \S 4]{freitag:1990}.  It is a phenomenon that only appears when $n>1$ that the functions obtained in this way are holomorphic.

\subsection{Poincare series}
Let $\mathbf{k} = (k_1, \cdots, k_n) \in \mathbb{Z}^n$ where $k_j>2$ for all $j= 1,\cdots,n$.  Let $\Gamma \subset \Gamma_F$ be a congruence subgroup.  For each $\nu \in \mathfrak{t}_\Gamma^\vee$ define the $\nu$'th \emph{Poincare series of weight $\mathbf{k}$} for $\Gamma$ by $$P_\nu(z) = \sum_{\gamma \in \Gamma_\infty \backslash \Gamma} \prod_{j=1}^n (c_\gamma^{(j)}z_j + d_\gamma^{(j)})^{-k_j}e^{2\pi i \langle \nu, \gamma z\rangle}\quad \text{where} \quad \Gamma_\infty = \left\{ \begin{pmatrix}
    1 & x \\
    0 & 1
\end{pmatrix}\mid x \in \mathfrak{t}_\Gamma \right\}.$$  One can show that the Poincare series is absolutely convergent and converges uniformly in compact subsets of $\mathbb{H}^n$ \cite[pg. 52]{garrett:1989}.  We check that $P_\nu$ is a Hilbert modular form of weight $\mathbf{k}$; for any $\delta \in \Gamma$, observe
\begin{align*}
& (P_\nu \vert_\mathbf{k} \delta)(z)\\
&= \prod_{j=1}^n (c^{(j)}_\delta z_j + d^{(j)}_\delta)^{-k_j} \sum_{\gamma \in \Gamma_\infty \backslash \Gamma} \prod_{j=1}^n \left( c^{(j)}_\gamma\left( \frac{a^{(j)}_\delta z_j +b^{(j)}_\delta}{c^{(j)}_\delta z_j + d^{(j)}_\delta}\right) + d^{(j)}_\gamma\right)^{-k_j}e^{2\pi i \langle\nu, \gamma z \rangle}\\
&= \sum_{\gamma \in \Gamma_\infty \backslash \Gamma} \prod_{j=1}^n (c^{(j)}_{\gamma \delta} z_j + d^{(j)}_{\gamma \delta})^{-k_j} e^{2\pi i \langle \nu , \gamma z \rangle} = P_\nu (z).
\end{align*}
Moreover, $P_\nu$ vanishes at each cusp of $\Gamma$ \cite[pg. 52]{garrett:1989}.  Thus $P_\nu \in S_\mathbf{k}(\Gamma)$ is a cusp form.  

It can be shown that the Poincare series span cusp forms.  Indeed, if $f\in S_\mathbf{k}(\Gamma)$ is an arbitrary cusp form with Fourier expansion $$f(z) = \sum_{\nu \in \mathfrak{t}^\vee} A(\nu) e^{2\pi i \langle \nu , z \rangle}$$ then $$\langle f, P_\nu \rangle = A(\nu) \cdot \operatorname{covol}(\mathfrak{t}_\Gamma) \prod_{j=1}^n (4\pi v_j)^{1-k_j} (k_j - 2)! \quad \cite[\text{pg. 53}]{garrett:1989}.$$ 
It follows that the only $f$ in the orthogonal complement of the subspace of $S_\mathbf{k}(\Gamma)$ spanned by the Poincare series are those $f$ with all Fourier coefficients zero, i.e. $f$ is the zero function.

\section{Hilbert modular forms (adelic formulation)}
\label{section:HMFadelic}

\subsection{Preliminaries}

In this section, we gather some preliminary material on valuations, places and adeles of a global field. For more details, see Cassels and Frohlich \cite{casselsfrohlich:1967}.

\subsubsection{Adeles}

We begin with the definition of \emph{restricted product} of topological spaces. Let $\lbrace X_{\alpha} \rbrace_{\alpha \in I}$ be a collection of topological spaces. Let $Y_{\alpha} \subset X_{\alpha}$ be open in $X_{\alpha}$ for all but finitely many $\alpha \in I$. Define
$$ X := \lbrace (x_{\alpha}) \ | \ x_{\alpha} \in X_{\alpha} \text{ for all } \alpha \in I, \text{ and } x_{\alpha} \in Y_{\alpha} \text{ for all but finitely many } \alpha \in I \rbrace. $$
We define a topology on $X$ by taking a basis of open sets to be of the form $\prod U_{\alpha}$ where $U_{\alpha} \subset X_{\alpha}$ is open for all $\alpha$ and $U_{\alpha} = Y_{\alpha}$ for almost all $\alpha$. The topological space $X$ is called the \emph{restricted product} of $\lbrace X_{\alpha} \rbrace$ with respect to $Y_{\alpha}$. 
We note the following properties of the restricted product:
\begin{enumerate}
    \item If $S$ is a finite subset of $\alpha$'s including those for which $Y_{\alpha}$ is not defined, let 
    \begin{equation*}
    \begin{split}
        X_S & := \lbrace (x_{\alpha}) \ | \ x_{\alpha} \in Y_{\alpha} \text{ for } \alpha \notin S \rbrace \\
    & \cong \prod_{\alpha \in S} X_{\alpha} \prod_{\alpha \notin S} Y_{\alpha}.
    \end{split}
    \end{equation*}
    $X_S$ is open in $X$ and the induced topology is same as the product topology. 
    \item If $X_{\alpha}$'s are locally compact, and $Y_{\alpha}$'s are compact, then $X$ is locally compact. 
    \item Suppose that $X_{\alpha}$'s are measure spaces with measure $\mu_{\alpha}$ satisfying $\mu_{\alpha}(Y_{\alpha}) = 1$ whenever $Y_{\alpha}$ is defined.  Then we define the measure on  $X$ by letting measurable sets to be of the form $\prod M_{\alpha}$, where $M_{\alpha} \subset X_{\alpha}$ has finite $\mu_{\alpha}$-measure and $M_{\alpha}=Y_{\alpha}$ for almost all $\alpha$, and 
    $$ \mu ( \prod M_{\alpha} ) := \prod \mu_{\alpha}(M_{\alpha}).$$
\end{enumerate}

Now, let $F$ be a global field, i.e., a finite extension of $\Q$ or a finite extension of the function field $\F_{q}(t)$. In the following, we will define a \emph{place} for $F$ assuming it is a finite extension of $\Q$. 

\begin{definition}
A \emph{valuation} on $F$ is a function $v : F \rightarrow \R_{\geq 0}$ such that:
\begin{enumerate}
    \item $v(\alpha) = 0$ if and only if $\alpha = 0$. 
    \item $v(\alpha \beta) = v(\alpha) v(\beta).$
    \item $\exists \ C$ such that $v(1+\alpha) \leq C$ whenever $v(\alpha) \leq 1$. 
\end{enumerate}
\end{definition}
If we can take $C=1$ in property 3, then $v$ is said to be non-archimedean. Otherwise, it is said to be archimedean. Two valuations $v_1$ and $v_2$ are said to be equivalent if and only if $v_1=v_2^c$ for some $c>0$. The valuation $v$ is said to be discrete if $\text{log} \circ v$ has discrete image in $\R$. 

Let $v$ be a nonarchimedean valuation on $F$. Let $\mathcal{O} := \lbrace \alpha \in F \ | \ v(\alpha) \leq 1 \rbrace$ be the valuation ring, $\mathcal{O}^{\times} := \lbrace \alpha \in F \ | \ v(\alpha)=1 \rbrace$ be the group of units, $\mathfrak{p} := \lbrace \alpha \in F \ | \ v(\alpha) < 1 \rbrace$ be the unique maximal ideal, and $\mathcal{O}/\mathfrak{p}$ be the residue field. Assume $v$ is discrete. Then, $\mathfrak{p}$ is a principal ideal $(\pi)$, where $\pi$ is a fixed uniformizer. (A necessary and sufficient condition for a nonarchimedean valuation to be discrete is that $\mathfrak{p}$ is a principal ideal.) The valuation $v$ is said to be normalized if $v(\pi) = 1/|\mathcal{O}_v/\mathfrak{p}|$. A valuation induces a metric, and hence a topology, on $F$. Let $F_v$ be the completion of the metric induced by $v$ on $F$. Then, the valuation $v$ on $F$ can be extended to $F_v$. The valuation ring $O_v$ of a complete field $F_v$ is compact. 

For an archimedean valuation $v$, if the corresponding completion of $F$ is $\R$, we say $v$ is normalized if it is equal to the ordinary absolute value on $\R$. If the completion is $\C$, we say it is normalized if it is equal to the square of the ordinary absolute value on $\C$.  

\begin{definition}
A \emph{place} of $F$ is an equivalence class of valuations on $F$. By a \emph{finite place} (resp. \emph{infinite place}) of $F$, we mean that the corresponding valuation is nonarchimedean (resp. archimedean). 
\end{definition}
Every place contains a unique normalized valuation. 

\begin{definition}
The \emph{finite adeles} of $F$, denoted by $\A_F^{\infty}$, is defined to be the restricted product of the completions $F_v$ with respect to $\mathcal{O}_v$ where $v$ is a finite place of $F$. The \emph{infinite adeles} of $F$, denoted by $\A_{F,\infty}$ is defined to be $F \otimes_{\Q} \R = \prod_{v | \infty} F_{v}$, where the product is over the infinite places of $F$. The \emph{adeles} of $F$, denoted by $\A_F$, is defined to be 
$$ \A_F := \A_F^{\infty} \times \A_{F,\infty}. $$
\end{definition}

$\A_F$ is a topological ring under component-wise addition and multiplication.
From property 2 of the the definition of restricted product and the fact that the valuation ring $\mathcal{O}_v$ is compact, it follows that $\A_F$ is locally compact. This implies, in particular, that the additive group of $\A_F$ has an invariant Haar measure. The field $F$ embeds into $\A_F$ under the diagonal map $x \mapsto (x,...,x,...)$ and the image is called the set of principal adeles. 

Let $K$ be a finite extension of $F$. Then, 
$$ \A_K \cong \A_F \otimes_F K$$
both algebraically and topologically. Lastly, we have the following important result.

\begin{theorem}[\cite{casselsfrohlich:1967}]
    The field $F$ is a discrete subgroup of $\A_F$ and the quotient $\A_F/F$ is compact. Further, the quotient measure induced on $\A_F/F$ is finite.
\end{theorem}

\begin{definition}
The units of the ring $\A_F$ are called the \emph{ideles} of $F$. They form a group under multiplication which is denoted by $\A_F^*$. The topology on $\A_F^*$ comes from embedding it in $\A_F \times \A_F$ by the map $x \mapsto (x,x^{-1})$. 
\end{definition}

When $F=\Q$, we write $\A$ in place of $\A_{\Q}$. Similar notation is employed for the finite and infinite parts.

\subsection{Definitions}

We primarily follow the exposition of Jarvis \cite{jarvis:1997} and Hida \cite{hida:1988}. Let $F$ be a fixed totally real field $F$ of degree $d$ over $\Q$. We will define the adelic version of the space of Hilbert modular forms for $F$ of weight $k = (k_1,...,k_d)$ and level $\mathfrak{n}$ (an ideal of $\mathcal{O}_F$). 

Let $I$ denote the set of embeddings of $F$ in $\R$.  Let $G$ be the linear algebraic group $\mathrm{Res}_{F/\Q}(\mathrm{GL}_2)$. This implies that if $A$ is a $\Q$-algebra, then the collection of $A$-points of $G$ is $G(A)=\mathrm{GL}_2(A \otimes_{\Q} F)$. Hence, $G(\A)=G(\A_{\Q})=\mathrm{GL}_2(\A_F)$ and similarly for the finite and infinite parts of the adele ring.

Let $t = (1,1,...,1) \in \Z^{I}$. We say $m \in \Z^I$ is \emph{parallel} if and only if $m \in \Z \cdot t$ and $m$ is parallel to $n$ if and only if $m-n \in \Z \cdot t$. We fix $k \in \Z^{I}_{\geq 1}$ such that 
$$ k_1 \equiv k_2 \equiv ... \equiv k_d \text{ (mod 2)}.$$
We define $v \in \Z^{I}_{\geq 0}$ such that $v_{\tau}=0$ for some $\tau \in I$ and $k+2v$ is parallel, and choose $w$ parallel to $v+k$. Such a $v$ is uniquely defined. We assume $w=v+k-t$ and let $\hat{w}=k-w$.

Let $f: G(\A) \rightarrow \C$ and $u=u^{\infty}u_{\infty} \in G(\A^{\infty}) \times G(\A_{\infty})_{+}$. Let $z_0 := (i,i,...,i) \in \mathbb{H}^d$. For $x = \left( \begin{bmatrix} a_{\tau} & b_{\tau} \\ c_{\tau} & d_{\tau} \end{bmatrix} \right)_{\tau \in I} \in G(\A_{\infty}) = \mathrm{GL}_2(\R)^d$, we define the automorphy factor 
$$ j(x,z) := (c_{\tau} z_{\tau} + d_{\tau} )_{\tau \in I}.$$
We define
\begin{equation}
(f|_{k,w} u)(x) := j(u_{\infty},z_0)^{-k} \mathrm{det}(u_{\infty})^w f(x u^{-1}).
\end{equation}

Let $U$ be an open compact subgroup of $G(\A^{\infty})$. We consider the following conditions on $f : G(\A) \rightarrow \C$.
\begin{enumerate}
    \item $f = f|_{k,w} u$ for all $u \in UC_{\infty,+}$, where $C_{\infty,+}$ is the stabilizer of $z_0=(i,i,...,i) \in \mathbb{H}^d$ in $G(\A_{\infty})_{+}$ which is identified with $\R^{\times} SO_2(\R)$.
    \item $f(ax)=f(x)$ for all $a \in G(\Q)=\mathrm{GL}_2(F)$.
    \item For each $z \in \mathbb{H}^d$, choose $u_{\infty} \in G_{\infty,+}$ such that $u_{\infty}(z_0)=z$. For $f$ satisfying the above conditions define $f_x : \mathbb{H}^d \rightarrow \C$ by
\begin{equation}
    f_x(z) := j(u_{\infty},z_0)^k \mathrm{det}(u_{\infty})^{-w} f(x u_{\infty})
\end{equation}
for all $x \in G(\A^{\infty})$. ($f_x$ is independent of the choice of $u_{\infty}$.) We require $f_x$ to be holomorphic with respect to each component $z_{\tau}$. 
\item We impose the condition:
$$ \int_{\A_F/F} f \left( \begin{bmatrix} 1 & a\\ 0 & 1 \end{bmatrix} x \right) dx = 0 $$
for all $x \in G(\A)$ and for each additive Haar measure on $\A_F/F$.
\item When $F=\Q$, we impose the extra condition that $| \text{Im}(z)^{k/2} f_x(z) |$ is uniformly bounded on $\mathbb{H}$ for all $x \in \mathrm{GL}_2(\A_{\Q})$. 
\end{enumerate}

We define the space of \emph{modular forms} $M_{k,w}(U)$ to be the space of functions $f: G(\A) \rightarrow \C$ satisfying conditions 1,2 and 3 (and 5 when $F=\Q$), and the space of \emph{cusp forms} $S_{k,w}(U)$ to be the space of functions $f$ in $M_{k,w}(U)$ which also satisfy condition 4. 

We define the following adelic analogues of the classical group $\Gamma_1(N)$. Let $U_0 = \prod_{\mathfrak{q}} \mathrm{GL}_2(\mathcal{O}_{F,\mathfrak{q}})$, where $\mathfrak{q}$ runs over the finite primes of $F$. Let $\mathfrak{n}$ be an ideal of $\mathcal{O}_F$. Define
$$ U_1(\mathfrak{n}) :=  \biggl\{ \begin{bmatrix} a & b \\ c & d \end{bmatrix} \in U_0 \ | \ c \in \mathfrak{n}, a-1 \in \mathfrak{n} \biggr\},$$
$$V_1(\mathfrak{n}) := \biggl\{ \begin{bmatrix} a & b \\ c & d \end{bmatrix} \in U_0 \ | \ c \in \mathfrak{n}, d-1 \in \mathfrak{n} \biggr\}. $$

We denote $S_{k,w}(\mathfrak{n}) = S_{k,w}(U_1(\mathfrak{n}))$ and $S^{*}_{k,w}(\mathfrak{n}) = S_{k,\hat{w}}(V_1(\mathfrak{n}))$. These are the \emph{cusp forms of level $\mathfrak{n}$}. 

We will now describe the relation of the above definition of Hilbert modular forms with the definition given in Section \ref{section:Hilbertmodularforms} . 

We define a congruence subgroup $\Gamma$ of $\mathrm{GL}_2(F) \cap G(\A^{\infty}) G(\A_{\infty})_{+}$ to be a subgroup which contains 
$$ \Gamma_{\mathfrak{n}} := \lbrace \gamma \in \mathrm{SL}_2(\mathcal{O}_F) \ | \ \gamma - I_2 \in \mathfrak{n} M_2(\mathcal{O}_F) \rbrace$$
for some integral ideal $\mathfrak{n}$ and such that $\Gamma / (\Gamma \cap F)$ is commensurable with $\mathrm{SL}_2(\mathcal{O}_F)/ \lbrace \pm 1 \rbrace$. For $f: \mathbb{H}^d \rightarrow \C$ and $\gamma \in G(\A_{\infty})$ we define
\begin{equation}
(f||_{k,w}\gamma)(z) := j(\gamma,z)^{-k} \mathrm{det}(\gamma)^{w} f(\gamma(z)).
\end{equation}

For $\Gamma$ a congruence subgroup as above, we consider the space of modular forms $M_{k,w}(\Gamma)$ to be the space of functions $f : \mathbb{H}^d \rightarrow \C$ such that the following conditions are satisfied:
\begin{enumerate}
    \item $f=f||_{k,w}\gamma$ for all $\gamma \in \Gamma$.
    \item $f$ is holomorphic with respect to each $z_{\tau}$. 
\end{enumerate}

If $f \in M_{k,w}(\Gamma)$, then $f$ has a Fourier expansion 
$$ f(z) = \sum_{\xi} a(\xi,f) e_F(\xi z),$$
where $a(\xi,f) \in \C$ are the Fourier coefficients and $e_F(\xi z) = \mathrm{exp}(2 \pi i \sum_{i=1}^d \xi^{\tau_i} z_i)$. Here $\xi$ runs over all the totally positive elements of a lattice in $F$ and zero. 
If $f$ vanishes at all the cusps of $\Gamma$, then $f$ is called a cusp form and the space of such forms is denoted by $S_{k,w}(\Gamma)$. (See Section \ref{section:Hilbertmodularforms}.)

Let $h$ be the class number of $F$. By the strong approximation theorem, for any finite set of places $S$, we can choose elements $t_1,t_2,...,t_h \in \A_F^{*}$ such that $t_i^{\infty} \in \prod_{\mathfrak{q}} \mathcal{O}_{F,\mathfrak{q}}$, where $\mathfrak{q}$ runs over the finite primes of $F$, and $t_1 \mathcal{O}_F, t_2 \mathcal{O}_F,...,t_h \mathcal{O}_F$ form a complete system of representatives of the ideal classes, and $(t_i)_S=1$. Here $t \mathcal{O}_F$ is the fractional ideal of $F$ corresponding to $t \in \A_F^*$. We will take $S$ to be a set of places which include all the infinite places of $F$ and finite places which divide the level $\mathfrak{n}$. 

Let 
$$ x_i = \begin{bmatrix} t_i & 0 \\ 0 & 1 \end{bmatrix}, \hspace{12mm} x_i^{-\iota} = \begin{bmatrix} 1 & 0 \\ 0 & t_i^{-1} \end{bmatrix}.$$
Let $E$ be the set of totally positive units of $\mathcal{O}_F$. We define 
$$ \Gamma_i(\mathfrak{n}) := x_i^{-\iota} E U_1(\mathfrak{n}) G(\A_{\infty})_{+} x_i^{\iota} \cap G(\Q) = x_i E V_1(\mathfrak{n}) G(\A_{\infty})_{+} x_i^{-1} \cap G(\Q).$$
Then, we have the following bijections:
\begin{equation}
\label{eq:cuspformbijection}
\begin{split}
    S_{k,w}(\mathfrak{n}) & \longrightarrow \bigoplus_{i=1}^h S_{k,w}(\Gamma_i(\mathfrak{n})), \\
    S^{*}_{k,w}(\mathfrak{n}) & \longrightarrow \bigoplus_{i=1}^h S_{k,\hat{w}}(\Gamma_i(\mathfrak{n})).
\end{split}
\end{equation}
In the second bijection, $f$ maps to $(f_{x_i})$ and the inverse map is $(f_1,f_2,...,f_h) \mapsto f$ where $f$ is given by 
\begin{equation}
f(\alpha x_i \gamma) = (f_i||_{k,w} \gamma_{\infty})(z_0)
\end{equation}
for all $\alpha \in G(\Q), \gamma \in U_1(\mathfrak{n})G(\A_{\infty})_{+}$. Here we use the decomposition  $$G(\A)= \bigsqcup_{i=1}^{h} G(\Q) x_i U_1(\mathfrak{n}) G(\A_{\infty})_{+}.$$

\subsection{Hecke theory, Fourier expansion and eigenfunctions}

We now define the Hecke algebras for $S_{k,w}(\mathfrak{n})$ and $S_{k,w}^*(\mathfrak{n})$. 

Let $U$ and $U'$ be open compact subgroups of $G(\A^{\infty})$. Let $x \in G(\A^{\infty})$. We define the operator
\begin{equation*}
    \begin{split}
    [UxU']: M_{k,w}(U) \rightarrow M_{k,w'}(U')  \\
    f \mapsto \sum_{i} f|_{k,w} x_i
    \end{split}
\end{equation*}
where we decompose $UxU'=\bigsqcup Ux_i$. 

Now, let $U = U_1(\mathfrak{n})$. For a prime $\mathfrak{q}$ of $F$, we define 
$$ T_{\mathfrak{q}} := \biggl [ U \begin{pmatrix} 1 & 0 \\ 0 & \pi_{\mathfrak{q}} \end{pmatrix} U \biggr ],$$
where $\pi_{\mathfrak{q}}$ is an element of  $\A^{\infty}$ which has 1 at every finite place, except at the place corresponding to $\mathfrak{q}$ where it is a uniformizer for $\mathfrak{q}$. For $\mathfrak{a}$, a fractional ideal of $F$ coprime to $\mathfrak{n}$, we define
$$ S_{\mathfrak{a}} := \biggl [ U \begin{pmatrix} \alpha & 0 \\ 0 & \alpha \end{pmatrix} U \biggr ],$$
where $\alpha = \prod_{\mathfrak{q}} \pi_{\mathfrak{q}}^{v_{\mathfrak{q}}(\mathfrak{a})}$.

The \emph{Hecke algebra} $\mathbb{T}_{k,w}(\mathfrak{n})$ is defined to be the $\Z$-algebra in $\mathrm{End}(S_{k,w}(\mathfrak{n}))$ generated by $T_{\mathfrak{q}}$ for all primes $\mathfrak{q}$ of $F$ and $S_{\mathfrak{a}}$ for all integral ideals $\mathfrak{a}$ of $F$ coprime to $\mathfrak{n}$. 

\begin{lemma}[\cite{jarvis:1997}]
\label{lemma:heckeisomorphism}
There is a natural isomorphism
\begin{equation*}
    \begin{split}
        S_{k,w}(U) \cong S_{k,\hat{w}}(U^{\iota}) \\
        f(x) \mapsto f^*(x)=f(x^{-\iota}),
    \end{split}
\end{equation*}
where $\iota$ denotes the main involution on $\text{GL}_2$, given by $x x^{\iota} = \mathrm{det}(x)$. We also have the relation
\begin{equation}
\label{eq:heckeisomorphism}
 (f|[UxU])^*=f^*|[U^{\iota}x^{-\iota}U^{\iota}]
\end{equation}
for all $x \in U$. 
\end{lemma}

Lemma \ref{lemma:heckeisomorphism} allows us to define the Hecke algebra $\mathbb{T}_{k,w}^*(\mathfrak{n})$ on $S_{k,w}^*(\mathfrak{n})$ as the $\Z$-algebra generated by $T_{\mathfrak{q}}$ and $S_{\mathfrak{a}}$ in $\text{End}(S_{k,w}^*(\mathfrak{n}))$ which is compatible with equation \ref{eq:heckeisomorphism}.

Lastly, we define the operator $< \mathfrak{a} >$ on $M_{k,w}^*(\mathfrak{n})$ using the identity 
$$ S_{\mathfrak{a}} = N (\mathfrak{a})^{[2w-k]} < \mathfrak{a} >,$$
where, if $r \in \Z \cdot t$, we write $r=[r] \cdot t$. 

We will now consider the notions of $q$-expansion, integrality and congruence of modular forms.

We identify the $q$-expansion of a modular form $f \in M_{k,w}^*(\mathfrak{n})$ with the $h$ Fourier expansions of the corresponding $f_i \in M_{k,\hat{w}}(\Gamma_i(\mathfrak{n}))$ as in \ref{eq:cuspformbijection}. So, the $q$-expansion will have $h$ constant terms and one term corresponding to all other totally positive ideals. 

Since $\Gamma_i(\mathfrak{n})$ contains elements of the form $\begin{bmatrix} 1 & b \\ 0 & 1 \end{bmatrix}$ for all $b \in t_i \mathcal{O}_F$, we see that the translation module for $\Gamma_i(\mathfrak{n})$ is $t_i \mathcal{O}_F$, embedded in $\mathbb{R}^n$ via the $n$ real embeddings of $F$.  Thus the dual lattice for this translation module is $t_i^{-1} \mathfrak{d}^{-1}$ embedded in $\mathbb{R}^n$ where $\mathfrak{d}$ is the \emph{different ideal} for $F$. So, $f_i$ has a Fourier expansion of the form 
\begin{equation} f_i(z) = \sum_{\substack{\xi \in t_i^{-1} \mathfrak{d}^{-1}\\ \xi\gg 0 }} a_i(\xi) e_F(\xi z).
\end{equation}

Let $\Phi(v)$ denote the subfield of $\overline{\Q}$ generated by $$ x^v := \prod_{\tau \in I} \tau(x)^{v_{\tau}}$$ for all $x \in F$.
Let $\mathcal{O}(v)$ be the ring of integers of $\Phi(v)$. The map $(-)^v: F^{\times} \rightarrow \Phi(v)^{\times}$ extends to a character on $\A_F^{\times}$ by continuity.

\begin{definition}
\label{definition:Hidacondition}
    Let $A$ be an $\mathcal{O}(v)$-algebra contained in $\C$. We say that $A$ satisfies the \emph{Hida condition} if for all $x \in \A_F^{\infty}$, $(x^v \mathcal{O}(v))A$ is generated by a single element in $A$.
    A number field $K$ containing $\Phi(v)$ is said to satisfy the Hida condition if $\mathcal{O}_K$ satisfies the \emph{Hida condition}.
\end{definition}
For $\mathfrak{q}=x \mathcal{O}_F$ a prime ideal in $\mathcal{O}_F$, we fix a choice of generator $\{ \mathfrak{q}^v \}$ of $x^v A$. If $\mathfrak{a}$ is any ideal, we define $\{ \mathfrak{a}^v \}$ by $\prod_{\mathfrak{q}} \{ \mathfrak{q}^v \}^{v_{\mathfrak{q}}(\mathfrak{a})}$. Then, $f \in S_{k,w}^*(\mathfrak{n})$ has a Fourier expansion
\begin{equation}
f \biggl( \begin{bmatrix} y & x \\ 0 & 1 \end{bmatrix} \biggr) = |y^{\infty}|_{\A}^{\epsilon} y_{\infty}^{k-w} \sum_{\xi \in F_{+}^{\times}} a(\xi y \mathfrak{d},f) \{ (\xi y \mathfrak{d})^v \} \xi^{-v} e_F(\sqrt{-1} \xi y_{\infty}) e_F(\xi x),
\end{equation}
where $a(\xi y \mathfrak{d},f)=0$ unless $\xi y \mathfrak{d}$ is an integral ideal, and $w=k+v-\epsilon \cdot t$ for any $\epsilon \in \Z$.

We define the space of \emph{$A$-integral (adelic) cusp forms} to be $$S_{k,w}^*(\mathfrak{n};A) := \{ f \in S_{k,w}^*(\mathfrak{n}) \ | \ a(\mathfrak{a},f) \in A \text{ for each integral ideal } \mathfrak{a} \}.$$

Then, $$f_i(z) = c_{v,i} \sum_{\xi \in (t_i^{-1} \mathfrak{d}^{-1})} a(\xi t_i \mathfrak{d},f) \{ \xi^v \} \xi^{-v} e_F(\xi z),$$ where $c_{v,i} = \mathrm{Nm}^{F}_{\Q} (t_i)^{-\epsilon} \{ (t_i \mathfrak{d})^v \}$. 

We also define the space $S_{k,\hat{w}}(\Gamma_i(\mathfrak{n});A)$ to be the space of modular forms $f \in S_{k,\hat{w}}(\Gamma_i(\mathfrak{n}))$ with Fourier coefficients in $A$. Hence, we see that
\begin{equation}
    S_{k,w}^*(\mathfrak{n};A) = \bigoplus_{i=1}^h c_{v,i} S_{k,\hat{w}} (\Gamma_i(\mathfrak{n});A)  
\end{equation} 

Similarly, for number fields $K$ containing $\Phi(v)$ and satisfying the Hida condition, we define the space of \emph{$K$-rational cusp forms} to be
$$S_{k,w}^*(\mathfrak{n};K) := \{ f \in S_{k,w}^*(\mathfrak{n}) \ | \ a(\mathfrak{a},f) \in K \text{ for each integral ideal } \mathfrak{a} \}.$$
We require this notion in Theorem \ref{theorem:mainresult}.

In \cite{hida:1988}, it is shown that if $A$ is an integrally closed domain containing $\mathcal{O}(v)$ satisfying the Hida condition and such that $A$ is finite flat over $\mathcal{O}(v)$ or $\Z_l$ for some prime $l \in \Z$, then $S_{k,w}^*(\mathfrak{n};A)$ is stable under the operator $< \mathfrak{a} >$ where $\mathfrak{a}$ is an integral ideal coprime to $\mathfrak{n}$. (This is true at least when $w=v+k-t$.) Similarly, $S_{k,w}^*(\mathfrak{n};K)$ is stable under the action of the Hecke operators $T_{\mathfrak{q}} \in \mathbb{T}_{k,w}^*(\mathfrak{n})$, at least when $w=v+k-t$. 

Let $K$ be a number field containing $\Phi(v)$ and satisfying the Hida condition and let $\mathcal{O}_K$ be its ring of integers. Let $\alpha$ be an integral ideal of $K$. We say that $f,g \in M_{k,w}^*(\mathfrak{n};K)$ are \emph{congruent modulo} $\alpha$ if the Fourier coefficients $a(\mathfrak{a},f) - a(\mathfrak{a},g) \in \alpha$ for all integral ideals $\mathfrak{a}$ of $K$. Similarly, $f$ is said to be a \emph{modulo} $\alpha$ \emph{eigenfunction} of the Hecke operator $T$ if there exists $a \in \mathcal{O}_K$ such that $f|T \equiv a \cdot f \text{ (mod } \alpha)$.

\subsection{$L$-functions of Hilbert modular forms}

We define $L$-functions of adelic Hilbert modular forms as in \cite{shimura:1978}.

Let $f\in S^{*}_{k,w}(\mathfrak{n})$ be an adelic Hilbert modular form corresponding to the $h$-tuple $(f_1,\cdots, f_h)$ of holomorphic Hilbert modular forms under the bijection in \ref{eq:cuspformbijection}. Here $f_i$ is a cusp form for $\Gamma_i(\mathfrak{n})$. 
We recall that $f_i$ has a Fourier expansion of the form 
$$ f_i(z) = \sum_{\substack{\xi \in t_i^{-1} \mathfrak{d}^{-1}\\ \xi\gg 0 }} a_i(\xi) e_F(\xi z).$$

Since $\Gamma_i(\mathfrak{n})$ contains elements of the form $\begin{bmatrix} \epsilon & 0 \\ 0 & 1 \end{bmatrix}$ for all totally positive units $\epsilon \in \mathcal{O}_F^*$, we have $$f_i(\epsilon^{(1)} z_1, \cdots \epsilon^{(n)} z_n) = \prod_{j=1}^n (\epsilon^{(j)})^{w_j} \cdot f_i(z) = \prod_{j=1}^n (\epsilon^{(j)})^{k_j/2} \cdot f_i(z).$$  The second equality above comes from the fact that $2w$ is parallel to $k$ and $\mathrm{Nm}^{F}_{\Q}(\epsilon)=1$. By taking the $\xi$'th Fourier coefficient of the above equality, we see
$$ a_i(\xi \epsilon^{-1})=  \prod_{j=1}^n (\epsilon^{(j)})^{k_j/2} \cdot a_i(\xi).$$  Thus the expression $a_i(\xi)\cdot \prod_{j=1}^n (\xi^{(j)})^{-k_j/2}$ for $0\ll \xi \in t_i^{-1}\mathfrak{d}^{-1}$ does not change when $\xi$ is replaced by $\epsilon \xi$ for any $0 \ll \epsilon \in \mathcal{O}_F^*$.

Now let $\mathfrak{m} \subset \mathcal{O}_F$ be an arbitrary ideal.  There exists a unique $\lambda \in\{1, \cdots, h\}$ such that $\mathfrak{m}$ and $t_\lambda^{-1} \mathcal{O}_F$ represent the same ideal class. If $\mathfrak{m}\neq 0$, it is possible to choose $0 \ll \xi\in t_\lambda^{-1}\mathfrak{d}$ such that $\mathfrak{m} = \xi t_\lambda^{-1}\mathcal{O}_F$.  Moreover, this choice of $0 \ll \xi\in t_\lambda^{-1}\mathfrak{d}$ is determined uniquely up to multiplication by totally positive units $\epsilon \in \mathcal{O}_F$.  Thus we may define $c(\mathfrak{m}, f)$ for any nonzero ideal $\mathfrak{m}\subset \mathcal{O}_F$ as $$c(\mathfrak{m}, f)  = a_\lambda(\xi) \cdot\prod_{j=1}^n (\xi^{(j)})^{-k_j/2}$$ where $\mathfrak{m} = \xi t_\lambda^{-1} \mathcal{O}_F$ and this is well defined.  For convenience, we define $c(\mathfrak{m},f) =0$ for $\mathfrak{m} =0$ and for non-integral fractional ideals as well. Let $C(\mathfrak{m},f) = N(\mathfrak{m})^{k_0/2}\cdot c(\mathfrak{m},f)$, where $k_0 = \max(k)$.

\begin{definition}
The \emph{L-function} attached to $f$ is given by 
\begin{equation}
\label{eq:dirichletseries}
D(s,f) = \sum_{\mathfrak{m} \subset \mathcal{O}_F} \frac{C(\mathfrak{m},f)}{N(\mathfrak{m})^s}
\end{equation}
where $\mathfrak{m}$ ranges over all ideals in $\mathcal{O}_F$. 
\end{definition}

As stated in \cite{shimura:1978}, the Dirichlet series in \ref{eq:dirichletseries} converges for sufficiently large $\Re(s)$ and has a meromorphic continuation to the whole plane $s\in\mathbb{C}$.  When $f$ is an eigenfunction for all the Hecke operators, and $c(\mathcal{O}_F, f) = 1$, we say $f$ is a \emph{normalized Hecke eigenform}.

\section{Galois representations}
\label{section:Galoisreps}

Galois representations are continuous representations of the absolute Galois group of a field $K$ in a vector space over $\C$ or $\Q_p$. In this section, we discuss basic notions related to Galois representations and explain what it means for Galois representations to be attached to Hilbert modular forms. 

\subsection{Preliminaries}

We first give an overview of some preliminary notions related to Galois representations. Our discussion follows the exposition in Chapter 9 of \emph{A First course in Modular forms} by Fred Diamond and Jerry Shurman \cite{diamondshurman:2005}.

\subsubsection{Frobenius elements}\label{section:Frobelt}

Let $F/\Q$ be a Galois number field. If $p \in \Z$ is any prime, then 
$$p \mathcal{O}_F = (\mathfrak{p}_1 ... \mathfrak{p}_g)^e, $$
for some prime ideals $\mathfrak{p}_1,...,\mathfrak{p}_g$ of $\mathcal{O}_F$. The numbers $e$ and $g$ are called the \emph{ramification index} and \emph{decomposition index} respectively. The \emph{residue fields} $\mathbb{F}_{\mathfrak{p}_i} := \mathcal{O}_F/\mathfrak{p}_i$ are all isomorphic to $\F_{p^f}$ for some $f$, called the \emph{inertial degree}.  Also, $efg=[F:\Q]$.

Suppose $\mathfrak{p}$ is a prime in $\mathcal{O}_F$ dividing $p$. The \emph{decomposition group} $D_{\mathfrak{p}}$ of $\mathfrak{p}$ is defined to be the subgroup of Gal$(F/\Q)$ which fixes $\mathfrak{p}$, i.e.,
$$ D_{\mathfrak{p}} := \lbrace \sigma \in \text{Gal}(F/\Q) \ | \ \mathfrak{p}^\sigma = \mathfrak{p} \rbrace.$$

The decomposition group has order $ef$. (Hence, it has index equal to $g$ in the Galois group). Since, $D_{\mathfrak{p}}$ fixes $\mathfrak{p}$, it acts on the residue field $F_\mathfrak{p}$. We get a group homomorphism $D_{\mathfrak{p}} \rightarrow \text{Gal}(\F_{\mathfrak{p}}/\F_p)$ which is, in fact, surjective.  The kernel of the map $D_{\mathfrak{p}} \rightarrow \text{Gal}(\F_{\mathfrak{p}}/\F_p)$ is defined to be the \emph{inertia group} of $\mathfrak{p}$, denoted by $I_{\mathfrak{p}}$.    
 
The order of $I_{\mathfrak{p}}$ is equal to $e$. Thus, $p$ is unramified if and only if $I_{\mathfrak{p}}$ is trivial. 

The quotient map $D_{\mathfrak{p}}/I_{\mathfrak{p}} \rightarrow \text{Gal}(\F_{\mathfrak{p}}/\F_p)$ is an isomorphism. We note that $\text{Gal}(\F_{\mathfrak{p}}/\F_p)$ is cyclic generated by the Frobenius automorphism $\sigma_p$, where $\sigma_p : x \mapsto x^p$. Then, we have the following definition.
\begin{definition}
A Frobenius element $\mathrm{Frob}_{\mathfrak{p}}$ of $\mathrm{Gal}(F/\Q)$ is any representative in $D_{\mathfrak{p}}$ whose coset in $D_{\mathfrak{p}}/I_{\mathfrak{p}}$ maps to $\sigma_p$. In other words, a Frobenius element of $\mathrm{Gal}(F/\Q)$ is any element $\mathrm{Frob}_{\mathfrak{p}} \in \text{Gal}(F/\Q)$ such that 
$$ x^{\mathrm{Frob}_{\mathfrak{p}}} \equiv x^p \ \ (\mathrm{mod} \ \ \mathfrak{p}), $$
for all $x \in \mathcal{O}_F$. 
\end{definition}
Note that, if $p$ is unramified, the Frobenius element $\text{Frob}_{\mathfrak{p}}$, associated to $\mathfrak{p}$, is unique. Moreover, since $F$ is Galois, if $\mathfrak{p}_1$ and $\mathfrak{p}_2$ are any two primes dividing $p$, then there is a $\sigma \in \text{Gal}(F/\Q)$ which sends $\mathfrak{p}_1$ to $\mathfrak{p}_2$. Hence, the subgroups $D_{\mathfrak{p}_1}$ and $D_{\mathfrak{p}_2}$ are conjugate in $\text{Gal}(F/\Q)$, which implies that the respective Frobenius elements $\text{Frob}_{\mathfrak{p}_1}$ and $\text{Frob}_{\mathfrak{p}_2}$ are also conjugate. We have the following weak version of the Chebotarev density theorem about Frobenius elements which generalizes the theorem of Dirichlet on primes in arithmetic progression. 

\begin{theorem}[\cite{diamondshurman:2005}]
\label{theorem:chebotarev}
Let $F/\Q$ be a Galois extension. Then, every element of $\mathrm{Gal}(F/\Q)$ takes the form $\mathrm{Frob}_{\mathfrak{p}}$ for infinitely many maximal ideals $\mathfrak{p}$ in $\mathcal{O}_F$.
\end{theorem}

\begin{remark}
Frobenius defined these Frobenius elements while investigating questions about splitting of primes in field extensions.    
\end{remark}

\subsubsection{Gal($\overline{K}/K$)}
Let $K$ be a field and write $\overline{K}$ to denote its separable closure. We review the profinite group structure of the absolute Galois group $G_{K} := \text{Gal}(\overline{K}/K$). 

Let $I$ be an indexing set with a partial order $\leq$. Suppose $\lbrace G_i \rbrace_{i \in I}$ is a collection of finite groups, equipped with the discrete topology, and $\phi_{ji}: G_j \rightarrow G_i$ for $i \leq j$ are homomorphisms satisfying $\phi_{ii} = \ $identity on $G_i$ and $\phi_{ki} \circ \phi_{jk} = \phi_{ji} , \ \ \forall i \leq k \leq j$. The collection $(\{G_i\}_{i\in I},\{\phi_{ji}\})$ is called an \emph{inverse system of groups}.  The \emph{inverse limit} $\varprojlim G_i$ of the inverse system $(\{G_i\}_{i\in I}, \{\phi_{ji}\})$ is defined to be the topological group
$$ \varprojlim G_i := \Big \lbrace (g_i)_{i \in I} \in \prod_{i \in I} G_i \ \ | \ \ \phi_{ji}(g_j) = g_i \ \ \forall i \leq j \Big \rbrace,$$
with underlying topology given by the relative product topology.  A topological group $G$ is called a \emph{profinite group} if it is isomorphic to the inverse limit of some inverse system of groups.  It is a consequence of the Tychonoff's theorem that all profinite groups are compact. 

The prototypical example of a profinite group is the group of $p$-adic integers $\Z_p$. It can be shown that $\Z_p$ is (isomorphic to) the inverse limit of $\Z/p^n\Z, $ for $n \geq 1$ with natural projections from $\Z/p^n\Z$ to $\Z/p^m\Z$ when $m \leq n$. 

Let $K/F$ be a Galois extension. We form an inverse system of groups Gal$(L/F)$ for finite, Galois subextensions $F \subset L \subset K$ with the natural projections $\text{Gal}(L_2/F) \rightarrow \text{Gal}(L_1/F)$ whenever $L_1 \subset L_2$. Then, $\text{Gal}(K/F) \cong \varprojlim \text{Gal}(L_i/F)$. The induced topology on $\text{Gal}(K/F)$ via this isomorphism is called the \emph{Krull topology}. The \emph{fundamental theorem of Galois theory} (for infinite extensions) states that there is an order-reversing bijection between closed subgroups $H$ of $\text{Gal}(K/F)$ and subextensions $F \subset L \subset K$. Under this bijection, a closed subgroup $H$  corresponds to the fixed field of $H$, and a subextension $L/F$ corresponds to the Galois group $\text{Gal}(K/L)$.  Moreover, open normal subgroups of $\text{Gal}(K/F)$ correspond, under this bijection, to finite, Galois subextensions. It is a fact that open normal subgroups form a neighborhood basis at the identity for a profinite group.  

We now carry over the discussion of Frobenius elements in Section \ref{section:Frobelt} to the setting of $\overline{\Q}/\Q$, where $\overline{\Q}$ is a fixed algebraic closure of $\Q$. Let $\overline{\Z}$ denote the ring of algebraic integers in $\overline{\Q}$. For $p \in \Z$, let $\mathfrak{p}$ be a maximal ideal of $\overline{\Z}$ lying over $p$. Let $\overline{\F_p}$ be an algebraic closure of $\F_p$, which we identify with the quotient $\overline{\Z}/\mathfrak{p}$. There is a natural reduction map $\overline{\Z} \rightarrow \overline{\mathbb{F}}_p$.  The \emph{decomposition group} $D_{\mathfrak{p}}$ is defined as
$$ D_{\mathfrak{p}} := \lbrace \sigma \in G_{\Q} \ | \  \mathfrak{p}^{\sigma} = \mathfrak{p} \rbrace.$$
The \emph{inertia group} $I_{\mathfrak{p}}$ is defined as the kernel of the group homomorphism $D_{\mathfrak{p}} \rightarrow G_{\overline{\mathbb{F}}_p}$.  This map is, in fact, surjective. Any pre-image $\text{Frob}_{\mathfrak{p}}$ of the Frobenius automorphism of $G_{\mathbb{F}_{p}}$ in $D_{\mathfrak{p}}$ is called an \emph{absolute Frobenius element} for the extension $\overline{\Q}/\Q$. A Frobenius element is only defined up to the inertia group. We also have that 
$$ \text{Frob}_{\mathfrak{p}_F} = \text{Frob}_{\mathfrak{p}} |_F,$$
where $F$ is any Galois number field and $\mathfrak{p}_F = \mathfrak{p} \cap F$. Lastly, if $\mathfrak{p}_1$ and $\mathfrak{p}_2$ are two maximal ideals lying over $p$, then the corresponding absolute Frobenius elements are conjugate in $G_{\Q}$. Then, by Theorem \ref{theorem:chebotarev} and the definition of Krull topology on $G_{\Q}$, we have the following theorem.

\begin{theorem}
\label{thereom:frobeniusdensity}
Let $S$ be a set consisting of all but finitely many primes in $\Z$. Then, the set of absolute Frobenius elements $\mathrm{Frob}_{\mathfrak{p}}$ corresponding to a family of maximal ideals $\mathfrak{p}$ lying over the primes in $S$ forms a dense subset of $G_{\Q}$. 
\end{theorem}

\subsection{Galois representations}

\subsubsection{$l$-adic Galois representations}

To understand more about the absolute Galois group $G_{\Q}$, we study its representations. One of the main motivations behind considering Galois representations is that they are a bridge between analytic objects such as modular forms and geometric objects such as elliptic curves. 

\begin{definition}
Let $K$ be either $\C$ or a finite extension of $\Q_{l}$ with the $l$-adic topology, where $l$ is a prime number. A Galois representation is a continuous group homomorphism $\rho : G_{\Q} \rightarrow \text{GL}(V)$, where $V$ is a finite-dimensional $K$-vector space. 
\end{definition}

When $K=\C$, the representation is generally called an \emph{Artin representation} and when, $K$ is a finite extension of $\Q_{l}$, it is called an \emph{$l$-adic Galois representation}. 

Examples of 1-dimensional Artin representations come from Dirichlet characters. Indeed, the composition of maps
$$ G_{\Q} \rightarrow \text{Gal}(\Q(\mu_N)/\Q) \rightarrow (\Z/N\Z)^* \rightarrow \C^*,$$
where the first map is restriction, second map is the isomorphism $(m_a: \zeta \mapsto \zeta^a) \mapsto a$ and the third is a primitive Dirichlet character $\chi$, gives a continuous homomorphism $\rho_{\chi}: G_{\Q} \rightarrow \C^*$. This has finite image. More generally, any continuous homomorphism from $G_{\Q}$ into $\C^*$ has finite image and factors through $\text{Gal}(F/\Q)$ for some Galois number field contained in a cyclotomic field.

To capture the infinite structure of $G_{\Q}$, one considers the \emph{$l$-adic cyclotomic character}. Let $\Q(\mu_{l^\infty})$ denote the field obtained by adjoining all the $l$-power roots of unity with $\Q$. Then, we have the following isomorphism
$$\text{Gal}(\Q(\mu_{l^\infty})/\Q) \cong \Z_{l}^*$$
by taking inverse limits on both sides of the isomorphism $\text{Gal}(\Q(\mu_{l^n})/\Q) \cong (\Z/l^n\Z)^*$. By pre-composing the isomorphism with the restriction $$G_{\Q} \rightarrow \text{Gal}(\Q(\mu_{l^\infty})/\Q),$$ we get a continuous homomorphism $\chi_l: G_{\Q} \rightarrow \Z_{l}^* \subset \Q_{l}^*$, which is called the $l$-adic cyclotomic character of $G_{\Q}$. As stated in \cite[eq. 9.13]{diamondshurman:2005}, we have  
$$ \chi_l (\text{Frob}_{\mathfrak{p}}) = p \text{ \ \ for\ \ } p \neq l,$$
and for any choice of $\mathfrak{p}$ and $\text{Frob}_{\mathfrak{p}}$, thus showing that the image of $\chi_l$ is infinite. 

Since $\text{Frob}_{\mathfrak{p}}$ is only defined up to inertia group, we need to impose the condition $I_{\mathfrak{p}} \subset \text{ker}(\rho)$ to get a well-defined notion of values of the representation $\rho$ at the Frobenius elements. This condition only depends on the underlying prime $p$ and the similarity class of $\rho$ because the inertia groups for different primes lying over $p$ are conjugate and $\text{ker}(\rho)$ is a normal subgroup. Furthermore, the characteristic polynomial of $\rho(\text{Frob}_{\mathfrak{p}})$ only depends on $p$  because the characteristic polynomial is invariant under conjugation. By Theorem \ref{thereom:frobeniusdensity}, the values $\rho(\text{Frob}_{\mathfrak{p}})$ determine $\rho$ completely. This discussion motivates the following definition.

\begin{definition}
Let $\rho$ be a Galois representation and $p \in \Z$ be a prime. Then $\rho$ is said to be unramified at $p$ if $I_{\mathfrak{p}} \subset \mathrm{ker}(\rho)$ for any maximal ideal $\mathfrak{p}$ of $\overline{\Z}$ lying over $p$. 
\end{definition}

One can see that $\rho$ being unramified at $p$ translates to the condition that $p$ is unramified in the field extension $K/\Q$ corresponding to the closed subgroup $\text{ker}(\rho)$ by Galois theory. 

We remark that the Artin representation $\rho_{\chi}$ which was defined using the primitive Dirichlet character $\chi$ is unramified at all $p$ not dividing $N$ and the $l$-adic cyclotomic character $\chi_l$ is unramified at all $p \neq l$. 

\subsection{$L$-functions of Galois representations}

In this section, we discuss the construction of $L$-functions attached to Galois representations of a number field $F$. In the case of $L$-functions attached to \emph{Artin representations} i.e. \emph{Artin $L$-functions}, our reference is \emph{Algebraic Number theory} by J{\"u}rgen Neukirch \cite{neukirch:1999}, and in the case of $l$-adic representations, we follow \emph{Abelian $l$-adic representations and Elliptic curves} by J. P. Serre \cite{serre:1989}. For simplicity, we only define the local factors of the $L$-functions at the non-archimedean places. 

\subsubsection{Artin $L$-functions}

Let $\rho: \text{Gal}(\overline{F}/F) \rightarrow \text{GL}(V) \cong \text{GL}_n(\C)$ be an Artin representation and let $\mathfrak{p}$ be a prime of $F$. Then, $\rho$ factors through some finite Galois extension $K/F$. Let $\mathfrak{P}$ be a prime of $K$ dividing $\mathfrak{p}$ and let $I_{\mathfrak{P}}$ be the inertia subgroup corresponding to $\mathfrak{P}$. We can choose a Frobenius element $\text{Frob}_{\mathfrak{P}} \in \text{Gal}(K/F)$. (It is uniquely defined if $\mathfrak{p}$ is unramified in $K$ i.e. if $I_{\mathfrak{P}}$  is trivial.)

We recall that $\rho$ is unramified  at the prime $\mathfrak{p}$ of $F$ if $I_{\mathfrak{P}} \subset \text{ker}(\rho)$. In this case, $\rho(\text{Frob}_{\mathfrak{P}})$ is well-defined. Let $V^{I_{\mathfrak{P}}}$ be the subspace of $V$ consisting of elements fixed under the action of $I_{\mathfrak{P}}$. It coincides with $V$ if $\rho$ is unramified at $\mathfrak{p}$.

We let $P_{\mathfrak{p},\rho}(T)$ denote the characteristic polynomial of $\rho^{I_{\mathfrak{P}}}(\text{Frob}_{\mathfrak{P}}) \in \C[T]$, where $\rho^{I_{\mathfrak{P}}}$ is the restriction of $\rho$ to the subspace $V^{I_{\mathfrak{P}}}$. The polynomial $P_{\mathfrak{p},\rho}(T)$ is well-defined because if $\rho$ is ramified at $\mathfrak{p}$, the inertia subgroup $I_{\mathfrak{P}}$ acts trivially on $V^{I_{\mathfrak{P}}}$, and hence $\rho^{I_{\mathfrak{P}}}(\text{Frob}_{\mathfrak{P}})$ is well-defined. Also, $P_{\mathfrak{p},\rho}(T)$ is independent of the choice of $\mathfrak{P}$ dividing $\mathfrak{p}$ because the different Frobenius elements corresponding to different primes $\mathfrak{P}$ dividing $\mathfrak{p}$ are conjugate to each other, and hence their characteristic polynomials are all the same.

\begin{definition}
The local $L$-factor corresponding to a prime $\mathfrak{p}$ of $F$ is defined to be
\begin{equation}
        L_{\mathfrak{p}}(\rho,T) 
        := \frac{1}{P_{\mathfrak{p},\rho}(T)} 
        = \frac{1}{\text{det}(I - \rho^{I_{\mathfrak{P}}}(\text{Frob}_{\mathfrak{P}}) T)}.
\end{equation}
\end{definition}

\begin{definition}
The $L$-function associated to $\rho$ is defined to be the Euler product 
\begin{equation}
    L(\rho,s) := \prod_{\mathfrak{p}} L_{\mathfrak{p}}(\rho, N(\mathfrak{p})^{-s})
    = \prod_{\mathfrak{p}} \frac{1}{\text{det}(I - \rho^{I_{\mathfrak{P}}}(\text{Frob}_{\mathfrak{P}}) N(\mathfrak{p})^{-s})}.
\end{equation}
\end{definition}

We also note that 
\begin{equation}
    \text{det}(I - \rho^{I_{\mathfrak{P}}}(\text{Frob}_{\mathfrak{P}}) N(\mathfrak{p})^{-s}) = \prod_{i=1}^{n} ( 1 - \epsilon_i N(\mathfrak{p})^{-s}),
\end{equation}
where the $\epsilon_i$'s are the eigenvalues of $\rho^{I_{\mathfrak{P}}}(\text{Frob}_{\mathfrak{P}})$. They are roots of unity because $\rho^{I_{\mathfrak{P}}}(\text{Frob}_{\mathfrak{P}})$ has finite order in $\text{GL}_n(\C)$. This defines $L(\rho,s)$ as an analytic function in the half-plane $\text{Re}(s)>1$. The \emph{Artin conjecture} states that $L(\rho,s)$ has an analytic continuation to the whole of complex plane for any non-trivial, irreducible $\rho$.

\subsubsection{$L$-functions associated to $l$-adic representations}

Let $F$ be a number field, and let $\Sigma_F$ denote the set of finite places of $F$. Let $\rho : \text{Gal}(\overline{F}/F) \rightarrow \text{GL}(V)$ be an $l$-adic representation of $F$ over a finite dimensional $\Q_l$-vector space $V$. Suppose $\mathfrak{p} \in \Sigma_F$ is unramified with respect to $\rho$. As before, set
$$ P_{\mathfrak{p},\rho}(T) := \text{det}(I -\rho(\text{Frob}_{\mathfrak{P}})T).$$

\begin{definition}
    The representation $\rho$ is said to be rational (resp. integral) if there exists a finite set $S \subset \Sigma_F$ such that
    \begin{enumerate}
        \item any element of $\Sigma_F \backslash S$ is unramified with respect to $\rho$, and
        \item if $\mathfrak{p} \notin S$, the coefficients of $P_{\mathfrak{p},\rho}(T)$ belong to $\Q$ (resp. to $\Z$). 
    \end{enumerate}
\end{definition}

\begin{definition}
    Let $\rho$ and $\rho'$ be $l$-adic and $l'$-adic representations of $F$ respectively. Assume $\rho$ and $\rho'$ are rational. Then $\rho$ and $\rho'$ are said to be \emph{compatible} if there exists a finite $S \subset \Sigma_K$ such that $\rho$ and $\rho'$ are unramified outside of $S$ and $P_{\mathfrak{p},\rho}(T) = P_{\mathfrak{p},\rho'}(T)$ for $\mathfrak{p} \in \Sigma_K \backslash S$. 
\end{definition}

\begin{definition}
    For each prime $l$, let $\rho_l$ be a rational $l$-adic representation of $F$. The system $\rho = \lbrace \rho_l \rbrace$ is said to be compatible if $\rho_l$ and $\rho_{l'}$ are compatible for any two primes $l$ and $l'$. The system $\rho=\lbrace \rho_l \rbrace$ is said to be \emph{strictly compatible} if there exists a finite $S \subset \Sigma_K$ such that the following conditions are satisfied. 
    \begin{enumerate}
        \item Let $S_l := \lbrace \mathfrak{p} \in \Sigma_K \ | \ \mathfrak{p} \text{ divides } l \rbrace$. Then, for all $\mathfrak{p} \notin S \cup S_l$, $\rho_l$ is unramified at $\mathfrak{p}$ and $P_{\mathfrak{p},\rho_l}(T)$ has rational coefficients.
        \item $P_{\mathfrak{p},\rho_l}(T) = P_{\mathfrak{p},\rho_{l'}}(T)$ for all $\mathfrak{p} \notin S \cup S_l \cup S_{l'}$. 
    \end{enumerate}
\end{definition}

The smallest such set $S$ is called the \emph{exceptional set} of the system $\rho = \lbrace \rho_l \rbrace$.

We will now define the $L$-function attached to a strictly compatible system of $l$-adic representations.  Let $\rho = \lbrace \rho_l \rbrace$ be a strictly compatible system of rational $l$-adic representations of a number field $F$ with  exceptional set $S \subset \Sigma_F$. For $\mathfrak{p} \notin S$, let $l$ be a prime $\neq p$ where $p$ is the rational prime below $\mathfrak{p}$. Since $\rho_l$ is unramified at $\mathfrak{p}$, and the polynomial $P_{\mathfrak{p},\rho_l}(T)$ has rational coefficients and is independent of the choice of $l$, we denote it simply by $P_{\mathfrak{p},\rho}(T)$. 

\begin{definition}
The $L$-function (at the finite unramified places) attached to $\rho$ is defined to be the Euler product
\begin{equation}
L(\rho,s) := \prod_{\mathfrak{p} \notin S} \frac{1}{P_{\mathfrak{p},\rho}(N(\mathfrak{p})^{-s})}.
\end{equation}
\end{definition}

As an example, take $\chi$ to be the system of $1$-dimensional $l$-adic cyclotomic characters $\chi_l$. Note, $\chi_l$ is ramified outside the set $\lbrace l \rbrace$, and $\chi_l(\text{Frob}_{\mathfrak{p}}) = p$ for any $l \neq p$ and for any choice of $\text{Frob}_{\mathfrak{p}}$ and prime $\mathfrak{p}$ dividing $p$. Hence, we see that $\chi = \lbrace \chi_l \rbrace$ is a strictly compatible system of $1$-dimensional integral $l$-adic representations of $\Q$ with the exceptional set $S$ as the empty set. The $L$-function attached to $\lbrace \chi_l \rbrace$ is 
$$ L(\chi,s) = \prod_{p} \frac{1}{(1-p^{1-s})} = \zeta(s-1).$$

One can similarly generalize the above definitions to the case of $K$-rational $\lambda$-adic representations, where $K$ is a number field and $\lambda$ is a prime of $K$, by replacing the instances of $\Q$ by $K$ and $l$ by $\lambda$.

\subsection{Attaching Galois representations to Hilbert modular forms}

We now state the main theorem which asserts the existence of a Galois representation associated to a Hilbert cuspidal eigenform.

\begin{theorem}[\cite{jarvis:1997}]
\label{theorem:mainresult}
    Let $f \in S_{k,w}^*(\mathfrak{n};K)$ be a $K$-rational Hilbert cuspidal eigenform, where $K$ is a number field containing $\Phi(v)$ and satisfying the Hida condition (see definition \ref{definition:Hidacondition}). Let $\theta_f(T)$ denote the eigenvalue of the Hecke operator $T$ acting on $f$. Let $\lambda$ be a prime of $\mathcal{O}_K$ lying above a prime $l \in \Z$ and let $\mathcal{O}_{K,\lambda}$ denote the completion of $\mathcal{O}_K$ at the prime $\lambda$. Then there exists a continuous representation
    $$ \rho_{\lambda} : \mathrm{Gal}(\overline{F}/F) \rightarrow \mathrm{GL}_2(\mathcal{O}_{K,\lambda})$$
    which is unramified outside $\mathfrak{n}l$ such that if $\mathfrak{q}$ is a prime of $\mathcal{O}_F$ not dividing $\mathfrak{n}l$, then
    \begin{equation}\label{eq:trRhoFrop}
    \mathrm{tr \ } \rho_{\lambda}(\mathrm{Frob}_{\mathfrak{q}}) = \theta_f(T_{\mathfrak{q}})
    \end{equation}
    and
    \begin{equation}\label{eq:delRhoFrob}
    \mathrm{det \ } \rho_{\lambda} (\mathrm{Frob }_{\mathfrak{q}}) = \theta_f(S_{\mathfrak{q}}) N(\mathfrak{q}).
    \end{equation}
\end{theorem}

Note that Equation \ref{eq:trRhoFrop} and Equation \ref{eq:delRhoFrob} imply that $\rho = \{ \rho_\lambda \}$ is a strictly compatible system of $K$-rational Galois representations with exceptional set contained in $\{\mathfrak{p} \in \Sigma_F \mid \mathfrak{p} \text{ divides } \mathfrak{n} \}$.

The theorem for $d=1$, i.e. $F=\Q$, was proved by Deligne \cite{deligne:1971} and Deligne-Serre \cite{deligne-serre:1974}. Carayol \cite{carayol:1986} proved it for odd $d$ whenever the weights are all at least 2,  and Taylor \cite{taylor:1989} and Blasius-Rogwaski \cite{blasiusrogawski:1989} proved it for even $d$ and all weights at least 2. When all the weights are equal to one, it was proved by Rogawski-Tunnell \cite{rogawski-tunnell:1983} and Ohta \cite{ohta:1984}. The last remaining case when some of the weights are equal to one was proved by Jarvis \cite{jarvis:1997}.

The method used by Taylor \cite{taylor:1989} and Jarvis \cite{jarvis:1997} is that of congruences. Taylor constructs congruences between the given form and a form of different suitable level, and then uses Wiles' method of `pseudo-representations' \cite{wiles:1988} to construct the desired representations.  Jarvis constructs congruences between forms of partial weight one and forms with higher weight thereby reducing the problem to the case where the theorem holds.

\bibliographystyle{amsalpha}
\bibliography{references}

\end{document}